\documentclass[12pt,a4paper]{article}
\usepackage[latin1]{inputenc}
\usepackage[T1]{fontenc}
\usepackage{amsmath,amssymb,amsthm,color}
\usepackage{graphicx}

\title{A generalization of the exponential sampling series and its approximation properties}
\author{
  Carlo Bardaro, Loris Faina and Ilaria Mantellini\\
 \small Department of Mathematics and Computer Sciences,\\
 \small University of Perugia,\\
 \small Via Vanvitelli, 1, I-06123 Perugia, Italy\\
 \small carlo.bardaro@unipg.it; loris.faina@unipg.it; mantell@dmi.unipg.it}
\date{}

\begin{document}
\maketitle
\noindent
{\small {\bf Abstract:} Here we introduce a generalization of the exponential sampling series of optical physics and establish pointwise and uniform convergence theorem, also in a quantitative form. Moreover we compare the error of approximation for Mellin band-limited functions using both classical and generalized exponential sampling series.}


%-----------------------------------------------------------------------------------------
%------------------------INTRODUZIONE----------------------------------------------------
\section{\bf Introduction}

The exponential sampling representation of a Mellin band-limited function, as a series in which the samples are exponentially spaced, is a powerful tool in finding solutions of certain inverse problems that have fundamental applications in optical physics phenomena, like the light scattering, Fraunhofer diffraction and radio astronomy (see \cite{OS}, \cite{BP}, \cite{GO}, \cite{CAS}). A mathematical rigorous version of the exponential sampling theorem for Mellin band-limited functions was firstly given in \cite{BJ0} and then considered in other aspects in \cite{BJ1}, \cite{BBM1}, \cite{BBM2}, especially with respect other basic formulae of Mellin transform analysis. Indeed, Mellin analysis is the more simple and suitable frame in which to put a rigorous theory of exponential sampling theory. It is important to remark that the exponential sampling formula can be formally viewed as the Mellin version of the classical Shannon sampling formula of Fourier analysis, using a suitable change of variables. However, this correspondance is in this case only apparent. As recently proved in \cite{BBMS} the notion of Mellin band-limited function is really different from the analogous one in Fourier analysis. Indeed, the two classes are disjoint, i.e. a function cannot be at the same time Mellin and Fourier band-limited. Also, the classical Paley-Wiener theorem of the Fourier analysis has a very different version in Mellin frame, which involves Riemann surfaces of the logarithm (see \cite{BBMS}). This fact motivates the study of an intrinsic theory of Mellin analysis, fully independent of the Fourier one, in which the results are deduced just using tools of Mellin transform theory and Mellin approximation theory. This project begun with the pioneering book \cite{MA}, developed in a systematic form  by P.L. Butzer and S. Jansche (see \cite{BJ2}, \cite{BJ1}) and then studied in a series of papers (see e.g. \cite{BM1}, \cite{BFM}, \cite{BBM3}, \cite{AV1}, \cite{AV2}). A well known development of the classical Shannon sampling theorem is given by the so-called generalized sampling theorem, in which the ``sinc'' kernel function is replaced by an arbitrary kernel function $\varphi$ satisfying suitable assumptions (see \cite{BS2}, \cite{BS1}, \cite{BBSV1}, \cite{BBSV2}). In this case the sampling theorem is given in an approximate sense, and the study of the error of approximation has a great importance in the applications, especially in prediction theory. Moreover, the use of a general kernel function enables one to approximate any function not necessarily (Fourier) band-limited.

In this paper, we introduce a  generalized version of the exponential sampling series acting on functions not necessarily Mellin band-limited, and we obtain pointwise and uniform convergence results, also in a quantitative form using suitable moduli of continuity. Also, we state an asymptotic formula of Voronovskaja type for locally regular functions $f.$ From a theoretical point of view, the approach used here is fully independent of Fourier analysis and uses only tools of Mellin analysis. We apply the theory to some particular cases, for which we compare the errors of approximation obtained for both the exponential sampling formula (truncation error) and its generalization (truncation and aliasing errors) applied to a Mellin band-limited function. Tables which show in details the behaviour of the two methods are given, for  two particular examples.

%-------------------------------------------------------------section 2--------------------------------------------------------------------------

\section{Basic notions and preliminary results}

Let $C(\mathbb{R}^+)$ be the space of all continuous and bounded functions defined on $\mathbb{R}^+.$ We say that a function $f \in C(\mathbb{R}^+)$ is "log-uniformly continuous" on $\mathbb{R}^+,$ if for any $\varepsilon >0$ there exists $\delta >0$ such that $|f(u) - f(v)| < \varepsilon$ whenever $|\log u - \log v| \leq \delta,$ for any $u,v \in \mathbb{R}^+.$ We denote by $\mathcal{C}(\mathbb{R}^+)$ the subspace of $C(\mathbb{R}^+)$ containing the log-uniformly continuous functions. Note that a log-uniformly continuous function is not necessarily uniformly continuous (in the usual sense), and conversely. Obviously the two notions are equivalent on compact intervals in $\mathbb{R}^+.$ 

We will say that a function $f$ is of class $C^{(n)}$ locally at the point $v \in \mathbb{R}^+$ if $f$ is $(n-1)$-times differentiable 
in a neighbourhood of $v$ and the derivative $f^{(n)}(v)$ exists.

For $1\leq p < +\infty,$ let $L^p= L^p(\mathbb{R}^+)$~ be the space of all the Lebesgue measurable and $p-$integrable complex-valued functions defined on $\mathbb{R}^+$ endowed with the usual norm $\|f\|_p.$ Analogous notations we give for functions 
defined on $\mathbb{R}.$

For $c \in \mathbb{R},$ let us consider the space
$$X_c = \{ f: \mathbb{R}^+\rightarrow \mathbb{C}: f(\cdot) (\cdot)^{c-1}\in L^1(\mathbb{R}^+) \}$$
endowed with the norm
$$ \| f\|_{X_c} = \|f(\cdot) (\cdot)^{c-1} \|_1 = \int_0^{+\infty} |f(u)|u^{c-1} du.$$

In an equivalent form, $X_c$ is the space of all functions $f$ such that $(\cdot)^c f(\cdot) \in L^1_\mu(\mathbb{R}^+),$ where $L^1_\mu= L^1_\mu(\mathbb{R}^+)$ denotes the Lebesgue space with respect 
to the (invariant) measure
$$\mu (A) = \int_A \frac{dt}{t},$$
for any measurable set $A \subset \mathbb{R}^+.$

The Mellin transform of a function $f\in X_c$ is defined by (see e.g. \cite{MA}, \cite{BJ2})
$$ [f]^{\wedge}_M (s) := \int_0^{+\infty} u^{s-1} f(u) du,~~(s=c+ it, t\in \mathbb{R}).$$
Basic properties of the Mellin transform are the following
$$[af(\cdot) + bg(\cdot)]^{\wedge}_M(s) = 
a [f]^{\wedge}_M(s) + b[g]^{\wedge}_M(s)~~~(f,g \in X_c,~a,b \in \mathbb{R})$$
$$|[f(\cdot)]^{\wedge}_M(s)| \leq \|f\|_{X_c}~~(s = c+it).$$
The inverse Mellin transform $M^{-1}_c[g]$ of the function $g \in L^1(\{c\} \times i\mathbb{R}),$ is defined by:
$$M^{-1}_c[g](x) \equiv M^{-1}_c[g(c+it)](x) := \frac{x^{-c}}{2 \pi}\int_{-\infty}^{+\infty} g(c+it) x^{-it}dt, ~~(x \in \mathbb{R^+}),$$
where by $L^p(\{c\} \times i\mathbb{R}),$ for $p \geq 1,$ we mean the space of all functions $g:\{c\} \times i\mathbb{R} \rightarrow \mathbb{C}$ with 
$g(c +i\cdot) \in L^p(\mathbb{R}).$

The Mellin translation operator $\tau_h^c$, for $h \in \mathbb{R}^+,~c \in \mathbb{R},$~$f: \mathbb{R}^+ \rightarrow \mathbb{C},$ is defined by
$$(\tau_h^c f)(x) := h^c f(hx)~~(x\in \mathbb{R}^+).$$
\noindent
Setting $\tau_h:= \tau^0_h,$ then  $(\tau_h^cf)(x) = h^c (\tau_hf)(x),$  
$\|\tau_h^c f\|_{X_c} = \|f\|_{X_c}.$

The Mellin convolution product of two functions $f,g: \mathbb{R}^+ \rightarrow \mathbb{C}$, denoted by $f\ast g,$  is defined by
$$ (f\ast g)(x) := \int_0^{+\infty} g(\frac{x}{u}) f(u) \frac{du}{u} = \int_0^{+\infty} (\tau^c_{1/u} f)(x)g(u)u^c \frac{du}{u}$$
in case the integral exists. For the properties of the Mellin convolution product see \cite{BJ2}.

The pointwise Mellin differential operator $\Theta_c,$ or the pointwise Mellin derivative $\Theta_cf$ of a function $f: \mathbb{R}^+ \rightarrow \mathbb{C}$ and $c \in \mathbb{R},$ is defined by (see \cite{BJ2})
$$
\Theta_cf(x) := x f'(x) + c f(x),~~x \in \mathbb{R}^+
$$
provided $f'$ exists a.e. on $\mathbb{R}^+.$ The Mellin differential operator of order $r \in \mathbb{N}$ is defined iteratively by
$$
\Theta^1_c := \Theta_c ,\quad\quad \Theta^r_c := \Theta_c (\Theta_c^{r-1}).
$$
For convenience set $\Theta^r:= \Theta^r_0$ for $c=0$ and $\Theta_c^0 := I,$ $I$ denoting the identity.
For instance, the first three Mellin derivatives are given by:
$$\Theta_cf(x) = xf'(x) + cf(x),$$
$$\Theta^2_cf(x) = x^2 f''(x) + (2c+1) xf'(x) + c^2f(x),$$
$$\Theta^3_cf(x) = x^3 f'''(x) + (3c+3)x^2f''(x) + (3c^2 + 3c +1)xf'(x) + c^3 f(x).$$
In general, we have (see \cite{BJ2})
$$\Theta^r_cf(x) = \sum_{k=0}^r S_c(r,k) x^k f^{(k)}(x),$$
where $S_c(r,k),$ $0\leq k\leq r,$ denote the generalized Stirling numbers of second kind, defined recursively by
$$S_c(r,0) := c^r,~ S_c(r,r) := 1,~ S_c(r+1,k) = S_c(r, k-1)+ (c+k) S_c(r,k).$$
 In particular for $c=0$
 $$\Theta^r f(x) = \sum_{k=0}^r S(r,k) x^k f^{(k)}(x)$$
 $S(r,k):= S_0(r,k)$ being the (classical) Stirling numbers of the second kind.
 
Also, in case $c=0$ we have the following Taylor formula with Mellin derivatives (see \cite{MA}, \cite{BM2}). 
For any $f \in C(\mathbb{R}^+)$ of class $C^{(n)}$ locally at the point $v,$ we have
$$f(tv) = f(v) + (\Theta f)(v) \log t + \frac{(\Theta^2f)(v)}{2!}\log^2t + \cdots \frac{(\Theta^nf)(v)}{n!}\log^nt 
+ h(t)\log^nt,$$
where $h:\mathbb{R}^+ \rightarrow \mathbb{R}$ is a bounded function such that $h(t) \rightarrow 0$ for $t \rightarrow 1.$

%-----------------------------section-------------------
\section{The generalized exponential sampling series}
 
Let $B^1_{c,T}$ denote the class of functions $f\in X_c,~f \in C(\mathbb{R}^+),$ $c \in \mathbb{R},$ which are Mellin band-limited 
in the interval $[-T,T],$ $T \in \mathbb{R}^+,$ thus for which $[f]^\wedge_M(c+it) = 0$ for all $|t| > T.$ 

A mathematical version of the exponential sampling theorem, introduced by the electrical engineers/physicists M.Bertero, E.R. Pike \cite{BP} and F. Gori \cite{GO} (see also \cite{OS}) and proved in \cite{BJ0}, reads as follows
\newtheorem{Theorem}{Theorem}
\begin{Theorem}[Exponential Sampling Formula]
If $f \in B^1_{c, \pi T}$ for some $c \in \mathbb{R},$ and $T>0,$ then the series 
$$x^c \sum_{k=-\infty}^{\infty} f(e^{k/T})\mbox{\rm lin}_{c/T}(e^{-k}x^T)$$
is uniformly convergent in $\mathbb{R}^+,$ and one has the representation
$$f(x) = \sum_{k=-\infty}^{\infty} f(e^{k/T})\mbox{\rm lin}_{c/T}(e^{-k}x^T) \equiv E^c_Tf(x)~~(x \in \mathbb{R}^+).$$
\end{Theorem}
\vskip0,4cm
The $\mbox{lin}_c-$function for $c\in \mathbb{R},$ $\mbox{lin}_c : \mathbb{R}^+ \rightarrow \mathbb{R},$ is defined, for $x \in \mathbb{R}^+ \setminus \{1\},$ by
\begin{eqnarray*} 
\mbox{lin}_c(x) = \frac{x^{-c}}{2 \pi i} \frac{x^{\pi i} - x^{-\pi i}}{\log x} = 
x^{-c} \mbox{sinc}(\log x) = \frac{x^{-c}}{2 \pi }\int_{- \pi}^\pi x^{-it}dt,
\end{eqnarray*}
with the continuous extension $\mbox{lin}_c(1) := 1.$

Here, the "sinc" function, as usual, is defined by
$$\mbox{sinc}(u) := \frac{\sin (\pi u)}{\pi u},~~u\neq 0, ~\mbox{sinc}(0) = 1.$$
\vskip0,3cm

The above basic theorem represents a Mellin version of the classical Shannon sampling formula, which involves Fourier band-limited functions. 
As proved recently in \cite{BBMS}, the class of Mellin band-limited functions and Fourier band-limited functions are disjoint, in the sense that a 
non-trivial function cannot be at the same time Mellin and Fourier band-limited. Moreover, using the Paley-Wiener theorem in Mellin setting 
(see \cite{BBMS}), a Mellin version of the classical Bernstein spaces must be defined throught the Riemann surfaces of the logarithm. This points out the importance to develop a theory of the exponential sampling fully independent of the Fourier analysis. 
This was done in recent papers (see \cite{BBM1}, \cite{BBM2}).

%-----------------------------subsection--------------------

\subsection{Definition of the series}

The function $\mbox{lin}_c \not \in X_{\overline{c}}$ for any $\overline{c}.$ This fact makes the theory difficult. Thus, we now introduce a new formula, 
which reconstructs a signal $f,$ not necessarily Mellin band-limited, in an approximate sense, using a general kernel function $\varphi \in C(\mathbb{R}^+)$ 
belonging to $X_c$ for some $c \in \mathbb{R}$ 

Let $\varphi: \mathbb{R}^+ \rightarrow \mathbb{R}$ be a continuous function such that the following assumptions are satisfied
\begin{enumerate}
\item[(i)] For every $u \in \mathbb{R}^+,$ 
$$\sum_{k=-\infty}^\infty \varphi(e^{-k}u) = 1,$$
and
$$M_0(\varphi):= \sup_{u \in \mathbb{R}^+}\sum_{k=-\infty}^\infty |\varphi(e^{-k}u)| < +\infty.$$
\item[(ii)] It holds
$$\lim_{r \rightarrow +\infty}\sum_{|k-\log u| > r}|\varphi(e^{-k}u)|=0,$$
uniformly with respect $u \in \mathbb{R}^+.$
\end{enumerate}
We denote by $\Phi$ the class of all the functions satisfying (i) and (ii). For any $\varphi \in \Phi$ we define the sampling operator
$$(S_w^\varphi f)(x) = \sum_{k=-\infty}^\infty f(e^{k/w})\varphi (e^{-k}x^w)\quad \quad (x \in \mathbb{R}^+, w >0),$$
for any function $f:\mathbb{R}^+ \rightarrow \mathbb{R}$ for which the series is absolutely convergent for any $x.$

In the next section we will prove some approximation properties for the above sampling operator.

%-------------------------------------------------subsection 2-------------------------------------------------

\subsection{Pointwise and uniform convergence}

Our first result is the following pointwise convergence result.

\begin{Theorem}\label{pointwise}
Let $\varphi \in \Phi$ and let $f$ be a bounded function. If $x \in \mathbb{R}^+$ is a continuity point for 
$f$ then
$$\lim_{w \rightarrow +\infty}(S_w^\varphi f)(x) = f(x).$$
\end{Theorem}
{\bf Proof}. Let $x \in \mathbb{R}^+$ be a continuity point of $f.$ For a fixed $\varepsilon >0$, let $\delta >0$ 
be such that $|f(x) - f(e^{k/w})| < \varepsilon$ whenever $|k-w \log x| \leq \delta w.$
By (i) we can write
\begin{eqnarray*}
&&|(S^\varphi_wf)(x) - f(x)| = \bigg|\sum_{k=-\infty}^\infty \varphi(e^{-k}x^w)[f(e^{k/w}) - f(x)]\bigg|\\
&\leq& \bigg(\sum_{|k-w\log x|\leq \delta w} + \sum_{|k-w\log x|> \delta w}\bigg) |\varphi(e^{-k}x^w)||f(e^{k/w}) - f(x)|\\
&\leq& M_0(\varphi) \varepsilon + 2\|f\|_\infty \sum_{|k-w\log x|> \delta w} |\varphi(e^{-k}x^w)|.
\end{eqnarray*}
The assertion follows now by (ii) letting $w \rightarrow + \infty.$ $\Box$
\vskip0,3cm
The next result states uniform convergence.
\begin{Theorem}\label{uniform}
Let $\varphi \in \Phi$ and $f \in \mathcal{C}(\mathbb{R^+}).$ Then
$$\lim_{w \rightarrow +\infty}\|S_w^\varphi f - f\|_\infty = 0.$$
\end{Theorem}
{\bf Proof}. The proof follows by the same arguments of Theorem \ref{pointwise}, taking into account that if 
$f \in \mathcal{C}(\mathbb{R^+})$ then we can choose $\delta >0$ such that for $|k-w \log x| < \delta w$ one has
$|f(x)- f(e^{k/w})| < \varepsilon$ uniformly with respect $x.$ $\Box$
\vskip0,4cm
\noindent
Now we state a quantitative approximation result for functions $f \in \mathcal{C}(\mathbb{R^+})$ in terms of the following 
modulus of continuity (see \cite{BM1}):
$$\omega (f, \delta):= \sup\{|f(s) - f(t)|: |\log s - \log t|\leq \delta\},\quad \quad \delta >0.$$
Note that $\omega$ satisfies all the classical properties of a modulus of continuity. In particular we will employ the following one:
$$\omega (f, \lambda \delta)\leq (\lambda +1)\omega(f, \delta),$$
for every $\delta, \lambda >0.$

For $j \in \mathbb{N}$ we define the (algebraic) moment of order $j$ of a kernel $\varphi \in \Phi$  
as 
$$m_j(\varphi,x):= \sum_{k=-\infty}^\infty \varphi(e^{-k}x)\log^j(e^kx^{-1})= \sum_{k=-\infty}^\infty \varphi(e^{-k}x)(k-\log x)^j,$$
and the absolute moments of arbitrary order $\alpha > 0$
$$M_\alpha(\varphi,x) := \sum_{k=-\infty}^\infty |\varphi(e^{-k}x)||k-\log x|^\alpha,$$
whenever the corresponding series are convergent. Finally, we put
$$M_\alpha(\varphi) := \sup_{x \in \mathbb{R}^+} M_\alpha(\varphi, x).$$

We have the following quantitative estimate
\begin{Theorem}\label{omega}
Let $\varphi \in \Phi$ be such that $M_1(\varphi) <+\infty,$ and let $f \in \mathcal{C}(\mathbb{R}^+).$ Then for every $\delta>0$ we have
$$|(S_w^\varphi f)(x) - f(x)| \leq M_0(\varphi) \omega(f, \delta) + \frac{\omega(f,\delta)}{w\delta}M_1(\varphi)\quad \quad (x \in \mathbb{R}^+).$$
\end{Theorem}
{\bf Proof}. We have
$$|(S_w^\varphi f)(x) - f(x)|\leq \sum_{k=-\infty}^\infty |\varphi (e^{-k}x^w)|\omega \bigg(f, \bigg|\frac{k}{w}- \log x\bigg|\bigg),$$
and, for any $\delta >0,$ using the properties of $\omega$ we can write
\begin{eqnarray*}
&&|(S_w^\varphi f)(x) - f(x)|\leq \sum_{k=-\infty}^\infty |\varphi (e^{-k}x^w)|\bigg(1+\frac{|(k/w) - \log x|}{\delta}\bigg) \omega(f, \delta)\\&\leq&
M_0(\varphi) \omega(f,\delta) + \frac{\omega(f,\delta)}{w\delta}M_1(\varphi). \Box
\end{eqnarray*}
As a consequence we obtain the following Corollary
\newtheorem{Corollary}{Corollary}
\begin{Corollary}
Under the assumptions of Theorem \ref{omega} we have
$$|(S_w^\varphi f)(x) - f(x)|\leq C \omega(f, \frac{1}{w}),$$
for an absolute constant $C$ depending only on $\varphi.$
\end{Corollary}
{\bf Proof}. For any fixed $w>0$ choose $\delta = 1/w$ and put $C=M_0(\varphi) + M_1(\varphi).$ $\Box$
\vskip0,4cm
Next, for functions $f$ which are of class $C^{(n)}$ locally at a point $x \in \mathbb{R}^+$ we will state an asymptotic 
formula which gives a precise order of pointwise approximation. In order to do that,
we use a stronger condition in place of (ii): there exists $n \in \mathbb{N}$ such that
\begin{enumerate}
\item[(iii)] For any $0\leq j\leq n,$ $m_j(\varphi) := m_j(\varphi, x)$ are independent of $x.$  
\item[(iv)] $M_n(\varphi) < +\infty$ and 
$$\lim_{r \rightarrow +\infty}\sum_{|k-\log u| >r}|\varphi(e^{-k}u)||k-\log u|^n = 0$$
uniformly with respect $u \in \mathbb{R}^+.$
\end{enumerate}
Note that if $M_n(\varphi) < +\infty$ then (ii) is satisfied. Indeed, we have
$$\sum_{|k-\log u| >r}|\varphi(e^{-k}u)| \leq \frac{1}{r^n}M_n(\varphi).$$

\begin{Theorem} \label{voronov}
Let $\varphi \in \Phi$ be a function satisfying (i), (iii) and (iv). If $f \in C(\mathbb{R}^+)$ is of class $C^{(n)}$ locally at the point 
$x\in \mathbb{R}^+$ then 
$$(S_w^\varphi f)(x) - f(x) = \sum_{j=1}^n \frac{(\Theta^jf)(x)}{j!}\frac{m_j(\varphi)}{w^j} + o(w^{-n}) \quad \quad (w \rightarrow +\infty).$$
\end{Theorem}
{\bf Proof}. Using the Taylor formula of order $n$ of $f$ we have
\begin{eqnarray*} 
&&|(S_w^\varphi f)(x) - f(x)|\\
&=& \bigg|\sum_{k=-\infty}^\infty \varphi (e^{-k}x^w)\bigg[\sum_{j=1}^n \frac{(\Theta^jf)(x)}{j!}\log^j\bigg(\frac{e^{k/w}}{x}\bigg) +
h\bigg(\frac{e^{k/w}}{x}\bigg) \log^n \bigg(\frac{e^{k/w}}{x}\bigg)\bigg]\bigg|
\end{eqnarray*}
Now, for a fixed $j$ one has
$$\sum_{k=-\infty}^\infty \varphi(e^{-k}x^w) \frac{(\Theta^jf)(x)}{j!}\log^j\bigg(\frac{e^{k/w}}{x}\bigg) = \frac{(\Theta^jf)(x)}{j!}\frac{m_j(\varphi)}{w^j},$$
thus we have only to estimate the remainder. Put
$$I:= \bigg|\sum_{k=-\infty}^\infty \varphi (e^{-k}x^w)h\bigg(\frac{e^{k/w}}{x}\bigg) \log^n \bigg(\frac{e^{k/w}}{x}\bigg)\bigg|.$$
For a fixed $\varepsilon >0$ let $\delta >0$ be such that $|h(e^{k/w}x^{-1})|< \varepsilon$ whenever $|k-w\log x| \leq \delta w.$ Then
\begin{eqnarray*}
&&w^n I < \varepsilon \sum_{|k-w\log x| \leq \delta w}|\varphi (e^{-k}x^w)| |k-w\log x|^n \\&+& \|h\|_\infty \sum_{|k-w\log x| \geq \delta w}|\varphi (e^{-k}x^w)| |k-w\log x|^n\\ &\leq& M_n(\varphi)\varepsilon + \|h\|_\infty\sum_{|k-w\log x| \geq \delta w}|\varphi (e^{-k}x^w)||k-w\log x|^n.
\end{eqnarray*}
Thus, using condition (iv), we obtain the assertion. $\Box$
\vskip0,4cm
As a Corollary, we immediately obtain the following pointwise approximation result:
\begin{Corollary} \label{voronov2}
Under the assumptions of Theorem \ref{voronov} we obtain
$$\lim_{w \rightarrow +\infty}w [(S_w^\varphi f)(x) - f(x)] = \Theta^1f(x)m_1(\varphi).$$
If moreover $m_j(\varphi) = 0,$ for $1\leq j \leq n-1,$  we have
$$\lim_{w \rightarrow +\infty}w^n [(S_w^\varphi f)(x) - f(x)] = \frac{\Theta^nf(x)}{n!}m_n(\varphi).$$
\end{Corollary}

%-----------------------------------examples-------------------------
\section{Examples}
In this section we give some examples of kernel $\varphi \in \Phi,$ satisfying the assumptions employed in the previous section. In what follows, for a given function $g:\mathbb{R}\rightarrow \mathbb{C},$ its Fourier transform, when it exists, will be defined as
$$\widehat{g}(v) := \int_{-\infty}^{+\infty} g(x) e^{-ivx}dx \quad \quad (v \in \mathbb{R}).$$
\subsection{Mellin splines}
We begin with an important class of functions with compact support, which represents the analogue  in the Mellin setting of the classical central B-splines. For every fixed $n \in \mathbb{N}$ we define 
$$B_n(t):= \frac{1}{(n-1)!}\sum_{j=0}^n (-1)^j\left(\begin{array}{c}  n\\ j\end{array} \right) \bigg(\frac{n}{2} + \log t - j\bigg)_{+}^{n+1}\qquad (t \in \mathbb{R^+}),$$
where, for every $r \in \mathbb{R},$ $r_+$ denotes the positive part of the number $r.$
More generally, for  $c \in \mathbb{R}$ one can consider the functions
$$B_{c,n}(t) = t^{-c}B_n(t)\quad \quad (t \in \mathbb{R}^+).$$
Since the functions $B_{c,n}$ have  compact support in $\mathbb{R}^+,$ and they are continuous, we have that $B_{c,n} \in X_c,$ for every $c \in \mathbb{R}.$ Thus we can compute the Mellin transform $[B_{c,n}]^{\wedge}_M,$ for $s=c+iv,$
$$[B_{c,n}]^{\wedge}_M(s) = \int_0^{+\infty} B_{c,n}(t) t^{c+iv-1}dt =
\int_0^{+\infty} B_{c,n}(t) e^{(c+iv)\log t}\frac{dt}{t}.$$
Making the substitution $\log t = z,$ and putting $\widetilde{B_n}(z):= B_n(e^z)$ we obtain, for every $c,$
$$[B_{c,n}]^{\wedge}_M(s) = \widehat{\widetilde{B_n}}(-v).$$
 Now, we have
$$\widetilde{B_n}(t) = \frac{1}{(n-1)!}\sum_{j=0}^n (-1)^j\left(\begin{array}{c}  n\\ j\end{array} \right) \bigg(\frac{n}{2} + t - j\bigg)_{+}^{n+1},$$
which is the classical central $B-$splines of order $n.$ Thus, 
$$[B_{c,n}]^{\wedge}_M(c+iv) = \bigg(\frac{\sin (v/2)}{v/2}\bigg)^n.$$

One can show that for $c=0,$ putting $B_{0,n}(t)=:B_n(t)$ for $t \in \mathbb{R}^+,$
$$\sum_{k=-\infty}^\infty B_n(e^kx) =1, \quad\mbox{for every}\quad x>0,$$
using the Mellin-Poisson summation formula (see \cite{BJ3}, \cite{BJ0}). In our present setting we have that the function
$$g(x):= \sum_{k=-\infty}^\infty B_n(e^kx)$$
is ``recurrent'' in the sense that $g(ex) = g(x),$ for every $x \in \mathbb{R}^+.$ Hence we can write the Mellin-Fourier series of $g$  as
$$g \sim \sum_{k=-\infty}^\infty g_{k,M}x^{-2\pi k i},$$
where, for $k \in \mathbb{Z},$
$$g_{k,M}:= \int_{1/e}^e g(x) x^{2\pi k i-1}dx$$
are the Mellin-Fourier coefficients of $g.$ Following the same proof of Theorem 7.2 in \cite{BJ3}, one has $g_{k,M}= [B_n]^{\wedge}_M(2k\pi i),$ $k \in \mathbb{Z}.$ This implies the Mellin-Poisson summation formula in the form
$$\sum_{k=-\infty}^\infty B_n(e^kx) = \sum_{k=-\infty}^\infty [B_n]^{\wedge}_M(2k\pi i) x^{-2k\pi i}.$$
Analogously we obtain corresponding results for the functions $B_{c,n},$ $c\neq 0,$ using the Mellin transform on the line $s= c+iv, v\in \mathbb{R}.$

Concerning the assumptions on the moments, since the functions $B_n$ have compact support, all the absolute moments $M_k(B_n),$ $k=0,1,\ldots, $ are finite numbers (indeed the series have a finite numbers of non-zero terms) and the values of the algebraic moments can be deduced again from the Mellin-Poisson summation formula, applied to the function $\psi_j(u) := (i\log u)^jB_n(u),$ $j=1,2,\ldots.$ As an example for $c=0$, by differentiating under the sign of the integral, it is not difficult to see that
$$\frac{d^j}{dt^j} [B_n]^\wedge_M(it) = [\psi_j]^\wedge_M(it),$$
and so the Mellin-Poisson formula for the function $\psi$ reads as
$$\sum_{k=-\infty}^\infty \psi(e^{-k}x) = \sum_{k=-\infty}^\infty \frac{d^j}{dt^j} [B_n]^\wedge_M(2k\pi i) x^{-2k\pi i}$$
which implies
$$i^j m_j(B_n)= \sum_{k=-\infty}^\infty \frac{d^j}{dt^j} [B_n]^\wedge_M(2k\pi i) x^{-2k\pi i}.$$
One can see that all algebraic moments $m_j(B_{n+1},x),$ are independent of $x$ for every $j=1,\ldots,n$ (see e.g. \cite{BM4}), so we can apply Theorem \ref{voronov} to the spline $B_{n+1}$ obtaining an asymptotic formula of order $n.$

A particular case is the second order Mellin spline ($n=2$) defined by (for $c=0$)
 \begin{eqnarray*}
 B_2(x) := (1-|\log x|)_{+}= \left\{\begin{array}{lll} 1-\log x, & \quad 1<x<e\\ 
 1+\log x, & \quad e^{-1}< x < 1\\ 0, & \quad \mbox{otherwise} \end{array} \right.
 \end{eqnarray*}
 In this case, since $m_1(B_2) = 0,$ the asymptotic formula reduces to (see Corollary \ref{voronov2}):
 $$ [(S_w^{B_2} f)(x) - f(x)] = o(w^{-1}) \quad \quad (w \rightarrow +\infty).$$
  
%-------------------------------------------------------------
   
\subsection{Mellin-Fejer kernel}

Another basic example is given by the family of Mellin-Fejer kernels, defined in the general case by
$$F^c_\rho (x):= \frac{x^{-c}}{2\pi}\rho ~ \mbox{sinc}^2(\frac{\rho}{\pi} \log\sqrt{x}),\quad \quad c \in \mathbb{R},\quad \rho >0,~  x \in \mathbb{R}^+,$$
with the continuous extension $F^c_\rho(1) = \rho/(2 \pi).$
These kernels are not of compact support, but $F^c_\rho \in X_c.$ Their Mellin transforms in $X_c$ are given by
\begin{eqnarray*}
[F^c_\rho]^\wedge_M(c+iv) = \left\{\begin{array}{ll} 1-\displaystyle\frac{|v|}{\rho}, & \quad |v|\leq \rho\\ 0, & \quad |v| >\rho. \end{array} \right. 
\end{eqnarray*}
For example, if we take $c=0$ and $\rho = 1,$ we have 
$$[F^0_1]^\wedge_M(iv) = (1-|v|)\chi_{[-1, 1]}(v),$$
where $\chi_A$ denotes the characteristic function of the set $A.$
Putting, as in the previous example, 
$$g(x)= \sum_{k=-\infty}^\infty F^0_1(e^kx)\quad \quad (x >0),$$
the function $g$ is recurrent in the sense that $g(ex) = g(x),$ for every $x>0,$ and again by the Mellin-Poisson summation formula we get
$$\sum_{k=-\infty}^\infty F^0_1(e^kx)= \sum_{k=-\infty}^\infty [F^0_1]^\wedge_M(2k\pi i)x^{-2k\pi i} = [F_1^0]^\wedge_M(0) = 1,$$
so that assumption  (i) is satisfied. As to (ii), we have only to show that
$$\lim_{r \rightarrow +\infty}\frac{1}{2\pi}\sum_{|k-\log u|> r}\mbox{sinc}^2(\log \sqrt{e^{-k}u}) = 0$$
uniformly with respect $u \in \mathbb{R}^+.$

Now, we can assume that $\log u \neq k$ for any integer $k.$
First, we have
$$\frac{1}{2\pi}\sum_{|k-\log u|> r}\mbox{sinc}^2(\log \sqrt{e^{-k}u}) \leq \frac{2}{\pi}\frac{1}{\sqrt{r}}\sum_{k=-\infty}^{\infty}\frac{\sin^2(\frac{1}{2}(\log u -k))}{|\log u -k|^{3/2}}.$$
Now, for such $u$ there exists an integer $\nu$ such that $\nu < \log u < \nu +1$ and we can write
$$\sum_{k=-\infty}^{\infty}\frac{\sin^2(\frac{1}{2}(\log u -k))}{|\log u -k|^{3/2}} \leq 
1+\sum_{k \in \mathbb{Z}, k \neq \nu, \nu+1}\frac{\sin^2(\frac{1}{2}(\log u -k))}{|\log u -k|^{3/2}}.$$
The last series is easily bounded by 
$$2\sum_{k=1}^\infty \frac{1}{k^{3/2}}$$
which is convergent. Thus we deduce (ii).

Therefore, one can apply Theorems \ref{pointwise} and \ref{uniform}.  However $M_1(\varphi) = +\infty,$ so we cannot apply Theorems \ref{omega} and \ref{voronov}. In the next example we modify the Mellin-Fejer kernel in a suitable way in order to have an asymptotic expression.

%--------------------------------------------------------

\subsection{The Mellin-Jackson kernel}

Let us consider the generalized Mellin-Jackson kernel, which is defined, for $c \in \mathbb{R},$ by (see \cite{BM4})
$$J_{\gamma, \beta}(x) := d_{\gamma, \beta}~ x^{-c}\mbox{sinc}^{2\beta}\bigg(\frac{\log x}{2\gamma \beta \pi}\bigg),$$
where $x \in \mathbb{R}^+,~\beta \in \mathbb{N},~\gamma \geq 1,~d_{\gamma, \beta}$ is a normalization constant, i.e.
$$d_{\gamma, \beta}^{-1}:=\int_0^{+\infty}\mbox{sinc}^{2\beta}\bigg(\frac{\log x}{2\gamma \beta \pi}\bigg)\frac{du}{u}.$$
It is easy to see that $J_{\gamma, \beta} \in X_c$  and $\|J_{\gamma, \beta}\|_{X_c} = 1.$ Moreover it is known that 
$[J_{\gamma, \beta}]^\wedge_M(c+iv) =0$ for $|v| \geq 1/\gamma,$ thus $J_{\gamma, \beta}$ is Mellin band-limited. 

Let us assume $c=0$ (the general case is treated in the same way). The Mellin-Poisson summation formula takes the form
$$\sum_{k=-\infty}^\infty J_{\gamma, \beta}(e^kx) = \sum_{k=-\infty}^\infty [J_{\gamma, \beta}]^\wedge_M(2k\pi i)x^{-2k\pi i} = [J_{\gamma, \beta}]^\wedge_M(0) = 1$$ 
thus assumption (i) is satisfied. Also, (ii) is satisfied, using similar arguments of the previous example. Concerning the moments, using again the Mellin-Poisson summation formula for the derivatives as in the first example, one has $m_1(J_{\gamma, \beta}) = 0,$ and, for $\beta >3/2,$ (see \cite{BM4})
$$m_2(J_{\gamma, \beta}) = d_{\gamma, \beta}\int_0^{+\infty} \mbox{sinc}^{2\beta}\bigg(\frac{\log x}{2\gamma \beta \pi}\bigg) \log^2x\frac{dx}{x}=: A_{\gamma, \beta} < +\infty.$$
Moreover, $M_2(J_{\gamma, \beta}) < +\infty.$ Therefore the assmptions of Theorems \ref{pointwise}, \ref{uniform}, \ref{omega}, \ref{omega2}, \ref{voronov} are satisfied, with $n=2.$ In particular, denoting by $(J_{\gamma, \beta}f)(x)$ the corresponding generalized sampling series, we obtain the following Voronovskaja formula, for $f \in C(\mathbb{R}^+)$ of class $C^{(2)}$ locally at the point $x \in \mathbb{R}^+:$
$$\lim_{w \rightarrow +\infty}w^2[(J_{\gamma, \beta}f)(x) - f(x)] = A_{\gamma, \beta}
\frac{\Theta^2f(x)}{2}.$$
%------------------------------------------------------------

\section{Some numerical evaluations}

In this section we will compare the approximation to a specific Mellin band-limited function $f$ of the classical exponential sampling series and the generalized ones.
\begin{enumerate}
\item For any $\rho >0$ let us consider the function $f(x) = F^c_\rho(x)$ for $x >0.$
As we have seen $F^c_\rho$ is Mellin band-limited to the interval $[-\rho, \rho].$
So, taking into account that now $T = \rho/\pi,$ the exponential sampling formula takes on the concrete form
\begin{eqnarray*}
F^c_\rho(x) &=& \sum_{k=-\infty}^\infty F^c_\rho(e^{k\pi/\rho})\mbox{lin}_{c\pi/\rho}(e^{-k}x^{\rho/\pi}) \nonumber\\
&=& \frac{2 \rho}{\pi^3 x^c}\sum_{k=-\infty}^\infty \frac{\sin^2(k \pi/2)}{k^2}\mbox{sinc}(\frac{\rho}{\pi}\log x -k).
\end{eqnarray*}
For the truncation error:
$$(T_{\rho, N}F^c_\rho)(x) := \frac{2 \rho}{\pi^3 x^c} \sum_{|k| \geq N+1}\frac{\sin^2(k \pi/2)}{k^2}\mbox{sinc}(\frac{\rho}{\pi}\log x -k),$$
we have the following pointwise estimate (see \cite{BBM1})
$$|F^c_\rho(x) - (S_NF^c_\rho)(x)| \leq |(T_{\rho, N}F^c_\rho)(x)| \leq \frac{4\rho}{\pi^4} N^{-1}~~~~(N \geq 2 \max\{x^{-c}, \frac{\rho}{\pi}|\log x|),$$
where $(S_NF^c_\rho)(x)$ denotes the N-th partial sum of the exponential sampling series of $F^c_\rho.$

Note that for the particular value of $c=0$ the pointwise estimate is of the order $\mathcal{O}(N^{-2})$ (see \cite[page 56]{BBM1}). Again, in certain situations, the estimate may be also better as we will show later.

Working now with a generalized exponential sampling formula, we have to take into account a truncation error:
$$(E_NF^c_{\rho})(x) = \sum_{|k| \geq N+1} F^c_\rho(e^{k/w})\varphi (e^{-k}x^w)\quad \quad (x \in \mathbb{R}^+, w >0),$$
and the aliasing error
$$(R_wF^c_\rho)(x):=|F^c_\rho(x) - (S_w^\varphi F^c_\rho)(x)|.$$
As an example, let $c=0$ and $\rho = \pi.$ Then 
$$f(x) = F^0_\pi(x) = \frac{1}{2}\mbox{sinc}^2(\frac{1}{2}\log x),$$
and the $N-$ partial sum of the exponential sampling series is given by
$$(S_NF^0_\pi)(x) = \frac{2}{\pi^2}\sum_{|k| \leq N}\frac{\sin^2(k\pi/2)}{k^2}\mbox{sinc}(\log x - k).$$
For the truncation error we have
$$|(T_{\pi, N}F^0_\pi)(x)| \leq \frac{4}{\pi^3}N^{-1}, \quad (N \geq 2\max\{1, |\log x|\}).$$
Now let us take the kernel $\varphi(x)= B_2(x) = (1-|\log x|)_{+}.$ In this case, the generalized exponential sampling series for the function $F^0_\pi$ is given by
\begin{eqnarray*}
(S^{B_2}_wF^0_\pi)(x) &=& \sum_{k=-\infty}^\infty F^0_\pi(e^{k/w}) B_2(e^{-k}x^w)\\&=& \frac{2w^2}{\pi^2}\sum_{k=-\infty}^\infty \frac{\sin^2(k\pi/2w)}{k^2}(1-|w\log x - k|)_+\\&=&
\frac{2w^2}{\pi^2}\sum_{k \in I_{w,x}} \frac{\sin^2(k\pi/2w)}{k^2}(1-|w\log x - k|)_+
\quad \quad (x \in \mathbb{R}^+, w >0),
\end{eqnarray*}
where $I_{w,x}:= \{k \in \mathbb{Z}: w\log x -1 < k <w\log x +1\},$ contains a finite number of integers corresponding to the non-zero terms. In particular, if $x$ is such that $1<x^w<e$ we have that the non-zero terms are for $k=0$ and $1.$ Then for a fixed $w>0$ and  $1 < x < e^{1/w}$ we have
$$(S^{B_2}F^0_\pi)(x) = \frac{2w^2}{\pi^2}\{\frac{\pi^2}{4w^2}(1-w\log x) + 
\sin^2(\pi/2w)(w \log x)\}.$$
If $e^{-1}<x^w<1,$ then the non-zero terms are obtained for $k=-1$ and $k=0,$ thus
$$(S^{B_2}F^0_\pi)(x) = \frac{2w^2}{\pi^2}\{\frac{\pi^2}{4w^2}(-w\log x) + 
\sin^2(\pi/2w)(1+w \log x)\}.$$
For $0<x^w<e^{-1},$ there exists an index $j \in \mathbb{N}$ such that $e^{-j-1}<x^w<e^{-j}$ and therefore the non zero terms are given by $k=-j-1$ and $k=-j$, thus 
$$(S^{B_2}F^0_\pi)(x) = \frac{2w^2}{\pi^2}\{\frac{\sin^2 ((j+1)\pi/2w)}{(j+1)^2}(-w\log x-j) + \frac{\sin^2(j\pi/2w)}{j^2}(1+j+w\log x)\}.$$
Analogously, for $x^w>e$ there exists $j \in \mathbb{N}$ such that $e^j < x^w <e^{j+1}$ and so the non zero terms are given in correspondance of $k=j$ and $k=j+1,$ obtaining
$$(S^{B_2}F^0_\pi)(x) = \frac{2w^2}{\pi^2}\{\frac{\sin^2 (j\pi/2w)}{j^2}(1-w\log x+j) + \frac{\sin^2((j+1)\pi/2w)}{(j+1)^2}(w\log x - j)\}.$$
%%%%%%%%%%%%%%%%%%%%%%%%% LORIS %%%%%%%%%%%%%%%
The Tables 1 and 2 show the behaviour of the two exponential sampling formulas for two values of $x$, one bigger and one smaller than $1.0$: in the first case ($\log x=2.7$), $x^w$ is definitely greater than $e$, therefore the only two terms of the series in $(S^{B_2}_w F^0_\pi)(e^{2.7})$ which are different from zero are those corresponding to $k=[2.7 w], [2.7 w]+1$;
in the second case ($\log x=-0.6$),  $x^w$ is definitely smaller than $1/e$,
therefore the only two terms of the series in $(S^{B_2}_w F^0_\pi)(e^{-0.6})$ which are different from zero are those corresponding to $k=-[0.6 w], -[0.6 w]-1$. 
\begin{center}
\small
\begin{tabular}{|r|c|}
\multicolumn{2}{c}{\bf Table 1.a}\\ \hline
$N$ & $S_N F^0_\pi(x)$ \\ \hline
   $20$    & $0.0220621711295$\\
   $40$    & $0.0220673420431$\\
   $160$   & $0.0220680655055$\\
   $640$   & $0.0220680767936$\\
   $2560$  & $0.0220680769700$\\
   $10240$ & $0.0220680769728$\\
   $20480$   & $0.02206807697284$ \\
   $1310720$ & $0.02206807697284$ \\ \hline
\end{tabular}
\quad
\begin{tabular}{|r|c|}
\multicolumn{2}{c}{\bf Table 1.b}\\ \hline
$w$ & $S_w^{B_1} F^0_\pi(x)$\\ \hline
$16$ & $0.02203184447881$\\
$32$ & $0.02205460079620$\\
$128$ &$0.02206723658951$\\
$512$ &$0.02206802441556$\\
$2048$&$0.02206807368853$\\
$8192$&$0.02206807676757$\\
$16384$ &   $0.02206807693864$\\
$1048576$ & $0.02206807697284$\\ \hline
\end{tabular}
\par\vspace{0.5cm}
\centerline{The behaviour of the exponential sampling formulas for $\log x=2.7.$}
\end{center}

\begin{center}
\small
\begin{tabular}{|r|c|}
\multicolumn{2}{c}{\bf Table 2.a}\\ \hline
$N$ & $S_N F^0_\pi(x)$ \\ \hline
   $160$   & $0.3684198616659$ \\
   $320$   & $0.3684198642867$ \\
   $640$   & $0.3684198646143$ \\
   $1280$  & $0.3684198646552$ \\
   $2560$  & $0.3684198646603$ \\
   $5120$  & $0.3684198646610$ \\
   $10240$ & $0.3684198646611$ \\ \hline
\end{tabular}
\quad
\begin{tabular}{|r|c|}
\multicolumn{2}{c}{\bf Table 2.b}\\ \hline
$w$ & $S_w^{B_1} F^0_\pi(x)$\\ \hline
$512$ &           $0.3687223730165$\\
$4096$ &          $0.3685332821791$\\
$32768$ &         $0.3684387681848$\\
$1048576$ &       $0.3684203077163$\\
$8388608$ &       $0.3684199385036$\\
$33554432$      & $0.3684198692762$\\
$2199023255552$ & $0.3684198646611$\\ \hline
\end{tabular}
\par\vspace{0.5cm}
\centerline{The behaviour of the exponential sampling formulas for $\log x=-0.6.$}
\end{center}
From the above tables one can deduce that the error of approximation $S_N F^0_\pi- F^0_\pi$ 
behaves as $N^{-3}$ at the considered points.
\vskip0,3cm

%%%%%%%%%%%%%%%%%%%%%%%%  END LORIS %%%%%%%%%%%
\item As a second example we consider the same function $f$ as before and the Mellin-Jackson kernel
$$\varphi (u) = J_{1,2}(u) = C \mbox{sinc}^4\bigg(\frac{\log u}{4\pi}\bigg), \quad u >0,$$
where $C$ is a normalization constant. 

The corresponding generalized exponential sampling series takes now the form:
$$(S^{J_{1,2}}_wf)(u) = \frac{C}{2}\sum_{k=-\infty}^{\infty}\mbox{sinc}^2\bigg(\frac{k}{2w}\bigg)\mbox{sinc}^4\bigg(\frac{w\log x -k}{4\pi}\bigg), \quad x>0.$$
For a given $w>0$ the pointwise truncation error $(T_{w,N}f)(x)$ is given by
$$(T_{w,N}f)(x) = \sum_{|k| \geq N+1}\mbox{sinc}^2\bigg(\frac{k}{2w}\bigg)\mbox{sinc}^4\bigg(\frac{w\log x -k}{4\pi}\bigg).$$
Employing similar reasonings as before, choosing $N \geq 2\max\{1, |w\log x|\},$
\begin{eqnarray*}
|(T_{w,N}f)(x)| \leq Rw^2\sum_{k \geq N+1}\frac{1}{k^6} = \mathcal{O}(N^{-5})\quad (N \rightarrow +\infty)
\end{eqnarray*}
where $R$ is a suitable constant and $\mathcal{O}$ depends on $w.$ Thus in evaluating the global error of approximation, one has to consider also the behaviour of $(T_{w,N}f)(x)$ for large values of $w.$
\end{enumerate}
\vskip0,3cm
\noindent
%%%%%%%%%%%%%%%%%%%%%%%%% LORIS %%%%%%%%%%%%%%%
The Tables 3 and 4 show the behaviour of the two exponential sampling formulas for the same two values of $x$ taken in the tables 1 and 2.
For the evaluation of $S_w^{J_{1,2}} F^0_\pi(x)$ we should first get a numerical
estimate for the normalization constant C:
$$
C^{-1} = \int_0^{+\infty} \mbox{sinc}^4(\frac{\log x}{4\pi})\frac{dx}{x}.
$$
An iterated trapezoidal formula gives for $C^{-1}$ the approximated
value $8.37757951289894$.

\begin{center}
\small
\begin{tabular}{|r|c|}
\multicolumn{2}{c}{\bf Table 3.a}\\ \hline
$N$ & $S_N F^0_\pi(x)$ \\ \hline
   $20$    & $0.0220621711295$ \\
   $40$    & $0.0220673420431$ \\
   $80$   &  $0.0220679852089$ \\
   $160$   & $0.0220680655055$ \\
   $320$   & $0.0220680755395$ \\
   $640$   & $0.0220680767936$ \\
   $1280$  & $0.0220680769504$ \\
   $2560$  & $0.0220680769700$ \\
   $5120$  & $0.0220680769725$ \\
   $10240$ & $0.0220680769728$ \\
   $20480$ & $0.0220680769728$ \\
   $40960$ & $0.0220680769728$ \\
   $81920$ & $0.0220680769728$ \\
   $163840$& $0.0220680769728$ \\
   $327680$& $0.0220680769728$ \\ \hline
\end{tabular}
\quad
\begin{tabular}{|r|c|}
\multicolumn{2}{c}{\bf Table 3.b}\\ \hline
$w=N$ & $S_w^{J_{1,2}} F^0_\pi(x)$\\ \hline
   $20$    & $0.0206901738326$ \\
   $40$    & $0.0216805928972$ \\
   $80$   &  $0.0219658314674$ \\
   $160$   & $0.0220418453359$ \\
   $320$   & $0.0220614368310$ \\
   $640$   & $0.0220664082084$ \\
   $1280$  & $0.0220676602406$ \\
   $2560$  & $0.0220679743972$ \\
   $5120$  & $0.0220680530799$ \\
   $10240$ & $0.0220680727685$ \\
   $20480$ & $0.0220680776929$ \\
   $40960$ & $0.0220680789243$ \\
   $81920$ & $0.0220680792322$ \\
   $163840$& $0.0220680793091$ \\
   $327680$& $0.0220680793284$ \\ \hline
\end{tabular}
\par\vspace{0.5cm}
\centerline{The behaviour of the exponential sampling formulas for $\log x=2.7.$}
\end{center}

\begin{center}
\small
\begin{tabular}{|r|c|}
\multicolumn{2}{c}{\bf Table 4.a}\\ \hline
$N$ & $S_N F^0_\pi(x)$ \\ \hline
   $20$    & $0.3684183377931$ \\
   $40$    & $0.3684196731679$ \\
   $80$   &  $0.3684198407044$ \\
   $160$   & $0.3684198616659$ \\
   $320$   & $0.3684198642867$ \\
   $640$   & $0.3684198646143$ \\
   $1280$  & $0.3684198646552$ \\
   $2560$  & $0.3684198646603$ \\
   $5120$  & $0.3684198646610$ \\
   $10240$ & $0.3684198646611$ \\
   $20480$ & $0.3684198646611$ \\
   $40960$ & $0.3684198646611$ \\ \hline
\end{tabular}
\quad
\begin{tabular}{|r|c|}
\multicolumn{2}{c}{\bf Table 4.b}\\ \hline
$w=N$ & $S_w^{J_{1,2}} F^0_\pi(x)$\\ \hline
   $20$    & $0.3638936971911$ \\
   $40$    & $0.3672463030803$ \\
   $80$   &  $0.3681212476795$ \\
   $160$   & $0.3683445829700$ \\
   $320$   & $0.3684009916855$ \\
   $640$   & $0.3684151657259$ \\
   $1280$  & $0.3684187182187$ \\
   $2560$  & $0.3684196074648$ \\
   $5120$  & $0.3684198299166$ \\
   $10240$ & $0.3684198855471$ \\
   $20480$ & $0.3684198994569$ \\
   $40960$ & $0.3684199029346$ \\ \hline
\end{tabular}
\par\vspace{0.5cm}
\centerline{The behaviour of the exponential sampling formulas for $\log x=-0.6.$}
\end{center}
\vskip0,2cm
Note that with the choice $w= N$ one can see that the error $S_w^{J_{1,2}}F^0_\pi - F^0_\pi$ 
behaves as $N^{-5}$ according to the theoretical estimate.
\vskip0,5cm
%%%%%%%%%%%%%%%%%%%%%%%%  END LORIS %%%%%%%%%%%
%{\bf Aknowledgments} The authors have been partialy supported by the Gruppo Nazionale Analisi Matematica, Probabilit%\'a e Applicazioni (GNAMPA) of the Istituto Nazionale di Alta Matematica (INdAM) and by the Department of Mathematics %and Computer Sciences of University of Perugia.

%--------------------------------------------references------------

\end{document}